\documentclass[12pt]{article}
\usepackage{amsmath,amsfonts,amssymb,amsthm,amscd}
\usepackage[]{hyperref}
\newtheorem{thm}{Theorem}[section]

\newtheorem{lem}{Lemma}[section]

\newcounter{alphthm}
\newtheorem{theo}[alphthm]{Theorem}

\newtheorem{pro}[alphthm]{Proposition}
\def\bd{\begin{definition}}
\def\ed{\end{definition}}
\def\bth{\begin{theorem}}
\def\eth{\end{theorem}}
\def\bpf{\begin{proof}}
\def\epf{\end{proof}}
\def\bco{\begin{corollary}}
\def\eco{\end{corollary}}
\def\ble{\begin{lemma}}
\def\ele{\end{lemma}}
\def\bpr{\begin{proposition}}
\def\epr{\end{proposition}}
\newcommand{\be}{\begin{equation}}
\newcommand{\ee}{\end{equation}}
\newcommand{\bes}{\begin{equation*}}
\newcommand{\ees}{\end{equation*}}
\newcommand{\br}{\begin{remark}}
\newcommand{\er}{\end{remark}}
\newcommand{\ben}{\begin{enumerate}}
\newcommand{\een}{\end{enumerate}}
\newcommand{\beq}{\begin{eqnarray}}
\newcommand{\eeq}{\end{eqnarray}}
\newcommand{\beqn}{\begin{eqnarray*}}
\newcommand{\eeqn}{\end{eqnarray*}}
\newcommand{\R}{I\!\! R}%For symbol of Euclidean space R
\def\nn{\nonumber}

\def\p{\partial}

\def\p{\partial}

\newtheorem{theorem}{Theorem}[section]
\newtheorem{lemma}[theorem]{Lemma}
\newtheorem{proposition}[theorem]{Proposition}
\newtheorem{corollary}[theorem]{Corollary}
\theoremstyle{definition}
\newtheorem{definition}[theorem]{Definition}

\theoremstyle{remark}
\newtheorem{remark}[theorem]{Remark}
\title{The Schwarzian derivative and conformal transformation on Finsler manifolds}
\author{\small B. BIDABAD\thanks{The corresponding author, bidabad@aut.ac.ir; behroz.bidabad@math.univ-toulouse.fr .}, \  F. SEDIGHI}
\date{ }
\begin{document}
\maketitle
\begin{abstract}
 Thurston, in 1986, discovered that the Schwarzian derivative has mysterious properties similar to the curvature on a manifold. After his work, there are several approaches to develop this notion on Riemannian manifolds.
Here, a tensor field is identified in the study of global conformal diffeomorphisms on Finsler manifolds as a natural generalization of the Schwarzian derivative. Then, a natural definition of a Mobius mapping on Finsler manifolds is given and its properties are studied. In particular, it is shown that Mobius mappings are mappings that preserve circles and vice versa. Therefore, if a forward geodesically complete Finsler manifold admits a Mobius mapping,  then the indicatrix is conformally diffeomorphic to the Euclidean sphere $ S^{n-1}$ in $ \mathbb{R}^n $. In addition, if a forward geodesically complete absolutely homogeneous Finsler manifold of scalar flag curvature  admits a non-trivial change of Mobius mapping, then it is a Riemannian manifold of constant sectional curvature.
\end{abstract}
\emph{Keywords;}  Finsler;    Schwarzian; Mobius; conformal; projective; concircular.
\section{Introduction}
Thurston in \cite{Th} claimes different curvatures on a manifold, measure the deviation of a curve or a manifold from being flat and Schwarzian derivative measures the deviation of a conformal map for being Mobius transformation. Therefore, by analogy, many qualitative constructions in differential geometry can be made.

The Schwarzian derivative appears also in many areas of complex analysis but it occurs first and foremost in the study of the Mobius mappings defined by $T(x)=\frac{ax+b}{cx+d}$, where $ad-bc\neq0$.
Historically, definition and properties of Schwarzian derivative were first introduced by Lagrange in \textquotedblleft Sur la construction des cartes g\'{e}ographiques\textquotedblright 1781.
Let $f$ be  a $C^{\infty}$ non-constant injective real function  on $\mathbb{R}$. The \textit{Schwarzian derivative}  is defined by
%\begin{equation}\label{sch}
% S(f(x))=\frac{f'''(x)}{f'(x)}-\frac{3}{2}(\frac{f''(x)}{f'(x)})^2 ,
%\end{equation}
\begin{equation}\label{Def;sch}
 S(f )=\frac{f''' }{f' }-\frac{3}{2}(\frac{f'' }{f' })^2 ,
\end{equation}
 where $ f',$  $ {} f'',$  $ {}  f''', $ are first, second, and third derivatives of $ f $ with respect to $ x\in\mathbb{R}  $.
  The expression \eqref{Def;sch} is ubiquitous and tends to appear in seemingly unrelated fields of mathematics: classical complex analysis, differential equations, one-dimensional dynamics, as well as, more recently, Teichm\"{u}ller theory, integrable systems and conformal field theory. It has been extended by several authors such as Osgood and Stowe \cite{OS}, Carne \cite{Ca}, etc.
A classical notation for $ S(f(x)) $ is $ \{f,x\} $, or $\{w,x\}$ if we write $w=f(x)$, and is due to Cayley in 1880.

 The Schwarzian derivative is occurred first as an operator  which is invariant under the all linear {Mobius mappings} in the sense that, $T$ is Mobius if and only if $ S(T\circ f)=S(f)$ or if and only if $S(T)=0$.

  Let $g$ be a real function for which the composition $ f\circ g $ is defined, we have
\begin{equation*}%\label{sch1}
 S(f\circ g)=S(g)+(S(f)\circ g)(g')^2.
\end{equation*}
 It follows that $ S(f)=S(g) $, implies $ f=T\circ g $ for some Mobius transformation $T$. To recall some geometric virtues of Mobius transformations, they map circles to circles and they are the only conformal maps of the sphere to itself. In particular, Mobius transformations are the only functions with vanishing Schwarzian.
B. Osgood and D. Stowe \cite{OS} introduced the Schwarzian tensor  $ B_g(\varphi) $ for two conformal Riemannian manifolds $ (M,g) $ and $ (M,\bar{g}) $ with $ \bar{g}=e^{2\varphi}g $ by
\begin{equation*}
B_g(\varphi)=Hess(\varphi)-d\varphi\otimes d\varphi-\frac{1}{n}(\Delta\varphi-\|grad\varphi\|^2)g.
\end{equation*}
%where $Hess(\varphi)$ is the Hessian of $\varphi$ and $grad\varphi$ is the gradient of $\varphi$.
   They also have defined the Schwarzian operator $f:(M,g)\to(M,\bar{g})$ by $S(f)=B_g(\varphi) $ where $ \varphi=log\|df\|$.
  Recently it's shown Schwarzian is also very useful in the study of Riemann-Finsler geometry. Meanwhile, one of the present authors in several joint works has studied the Schwarzian derivative for projective transformations in Finsler geometry and obtained a short proof for some known results, see for instance  \cite{BiS,SB}.

  In the present work, we identify a tensor which arises in the study of conformal change of metrics on a Finsler manifold as a natural generalization of the Schwarzian derivative.  We use a certain conformal parameter to define the Schwarzian tensor which follow, fairly directly, from the definition, and which have corresponding formulations in the classical setting.
  Meanwhile, we obtain the following theorem;
\begin{thm}\label{th;group}
 Let $(M,F)$ be a Finsler manifold. The set of Mobius transformations forms a subgroup of the conformal group   and contains the homotheties group of $(M,F)$.
\end{thm}
Among the others, after a joint work with Z. Shen \cite{BS} it is shown that the Mobius transformations on Finsler geometry are equivalent to the circle preserving or concircular transformations.
 \begin{thm}\label{th;Con=Mobgroup}
The group of conformal transformations of $(M,F)$ coincides with its Mobius group if and only if it maps all geodesic circles to geodesic circles.
\end{thm}
These theorems imply the following rigidity theorems.
\begin{thm}\label{Th;MobiusComplet}
If a forward geodesically complete Finsler manifold $(M,F)$ admits a Mobius mapping, then the indicatrix is conformally diffeomorphic to the Euclidean sphere $ S^{n-1} $ in $ \mathbb{R}^n $.
\end{thm}
\begin{thm}\label{Th;MobiusScalar}
Let $(M,F)$ be a forward geodesically complete absolutely homogeneous Finsler manifold of scalar flag curvature. If $(M,F)$ admits a nontrivial Mobius mapping, then it is a Riemannian manifold of constant sectional curvature.
\end{thm}
%In the case of Einstein Randers manifolds we have the following application of Schwarzian derivative.

\begin{thm}\label{Th;Randers}
Let $ (M,F) $ be a compact boundaryless Einstein Randers manifold with constant Ricci scalar and the projective parameter  $p$.\\
\textbullet If $ S(p)=0 $, then $ (M,F) $ is Berwaldian.\\
\textbullet If $ S(p)<0 $, then $ (M,F) $ is Riemannian.
\end{thm}
\section{Preliminaries and notations}
\subsection{Finsler structure}
Let $M$ be an $n$-dimensional connected smooth manifold. We denote by  $TM$ the tangent bundle and $\pi:TM_0\to M$, the fiber bundle of non-zero tangent vectors. A \textit{Finsler structure} on $M$ is a function $F:TM \to [0,\infty )$, with the following properties: $F$ is  $ C^\infty$ on $TM_0$; $F$ is positively homogeneous of degree one in $y$, that is $F(x,\lambda y)=\lambda F(x,y) $ for all positive $ \lambda $; The Hessian matrix of $ F^2 $, defined by  $ (g_{ij})=(1/2[\frac{\partial^2}{\partial y^i \partial y^j}F^2]) $, is positive definite on $ TM_0 $.
 A \textit{Finsler manifold} is a pair consisting of a differentiable manifold $ M $ and a Finsler structure $ F $ on $ M $ denoted here by $ (M,F) $.
   The hyper-surface $S\subset T_{x_0}M$ defined by $S:=\{y\in T_{x_0}M:F(x_{0},y)=1\}$  is called an \emph{indicatrix} in $x_0\in M$.
  Every Finsler structure $F$ induces a spray $ G=y^i\frac{\partial}{\partial x^i}-G^i(x,y)\frac{\partial}{\partial y^i},$ {} on $ TM $, where $ G^i(x,y)=\frac{1}{4}g^{il}\{[F^2]_{x^ky^l}y^k-[F^2]_{x^l}\} $, $G$ is a globally defined vector field on $ TM $. Differential equation of a geodesic in local coordinates is given by $ \frac{d^2x^i}{ds^2}+2G^i(x(s),\frac{dx}{ds})=0 $,  where $ s(t)=\int_{t_0}^{t}F(\gamma,\frac{d\gamma}{dr})dr , $ {} is the arc length parameter. One can observe that the pair $ \{\frac{\delta}{\delta x^i},\frac{\partial}{\partial y^i}\} $ forms a horizontal and vertical frame for $ TTM $, where $\frac{\delta}{\delta x^i}=\frac{\partial}{\partial x^i}-G^j_i\frac{\partial}{\partial y^j} $, {} $ G^j_i=\frac{\partial G^j}{\partial y^i} $ and $ 2G^i=\gamma^i_{jk}y^jy^k $, where
\be\label{Eq;ChristRie}
\gamma^i_{jk}=\frac{1}{2}g^{ih}(\frac{\partial g_{hk}}{\partial x^j}+\frac{\partial g_{hj}}{\partial x^k}-\frac{\partial g_{jk}}{\partial x^h}),
\ee
 are formal \textit{Christoffel symbols} of the second kind.
A Finsler structure  $F$ is called \emph{forward (resp. backward) geodesically complete}, if every  geodesic on an open interval $(a, b)$ can be extended to a geodesic on $(a,\infty)$ (resp. $(-\infty, b)$).
$F$ is said to be \emph{complete} if it is forward and backward complete.

 Let $f:M\to\mathbb{R}$ be a smooth function on an $n$-dimensional $(n\geqslant2)$ Finsler manifold $(M,F)$. At a point $p\in M$,   the \textit{gradient} vector field of $f$,  $\nabla f(p)=grad\,f(p)\in \pi^*TM$, is defined by $g_{grad\,f(p)}(\nu,grad\,f(p))=df_p(\nu),$ $ \forall \nu\in T_pM,$
where $df:=\frac{\partial f}{\partial x^i}dx^i$ is the differential of $f$.
In terms of a local coordinate system, we have $grad\, f:=f^i(x)\frac{\partial}{\partial x^i}\in \pi^*TM$, where $f^i(x)=g^{ij}(x,grad\, f(x))\frac{\partial f}{\partial x^j}$.
\subsection{Cartan connection and Koszul formula}
Here, a brief global approach to the Cartan connection is recalled for our further setting.
 Any point of  $ TM_0 $ will be denoted by $ z=(x,y) $ where $ x=\pi z \in M $ and $ y \in T_{\pi z}M $.
  By $ TTM_0 $ we denote the tangent bundle of $ TM_0 $ and  by $ \pi^*TM $ the pull back bundle of $ \pi $. Consider the canonical linear mapping $ \varrho:TTM_0 \to \pi^*TM $, where $ \varrho=\pi_*$  and $ \varrho \hat{X}=X$ for all $ \hat{X}\in \Gamma(TM_0)$. Locally we have
   $\varrho_z(\frac{\delta}{\delta x^i})_z=(\frac{\partial}{\partial x^i})_x $ and $ \varrho_z(\frac{\partial}{\partial y^i})_z=0 $. Let $ V_zTM $ be the set of vertical vectors at $ z \in TM_0 $, that is, the set of all vectors which are tangent to the fiber through $ z $. Equivalently, $ V_zTM=ker\pi_* $ where $ \pi_*:TTM_0 \to TM $ is the linear tangent mapping. Let $\nabla  $ be a linear connection on $ \pi^*TM $ the sections of pull back bundle $ \pi^*TM $,
   \begin{align*}\label{Def;CartanConn}
   \nabla:T_zTM_0\times\Gamma(\pi^*TM)&\to \Gamma(\pi^*TM), \\
  (\hat{X},v)&\mapsto\nabla_{\hat{X}}v,
    \end{align*}
      provided that there is a linear mapping $ \mu:TTM_0 \to \pi^*TM $, defined by $ \mu(\hat{X})=\nabla_{\hat{X}} \upsilon $, where $ \hat{X}\in TTM_0 $ and $ \upsilon $ is the canonical section of $ \pi^*TM $. The connection $ \nabla $ is said to be \textit{regular}, if $ \mu $ defines an isomorphism between $ VTM_0 $ and $\pi^*TM $. In this case, there is a horizontal distribution $ HTM $ such that we have the \textit{Whitney sum} $ TTM_0=HTM\oplus VTM $. The linear connection on $ \pi^*TM $ is said to be a \textit{Finsler connection}, if it is regular. It can be shown that the sets $ \{\frac{\delta}{\delta x^j}\} $ and $ \{\frac{\partial}{\partial y^j}\} $, form a local frame field for the horizontal and vertical subspaces and the dual frame $\lbrace dx^{i}\rbrace$ and $\lbrace\delta y^{i}\rbrace$ respectively. This decomposition permits to write a vector field $ \hat{X} \in TTM_0 $ into the horizontal and vertical form $ \hat{X}=H\hat{X}+V\hat{X} $, uniquely.

The \textit{torsion tensor} of the Finsler connection $ \nabla $ is defined by
\bes%\label{tor}
   \tau(\hat{X},\hat{Y})=\nabla_{\hat{X}}Y-\nabla_{\hat{Y}}X-\varrho[\hat{X},\hat{Y}].
 \ees
     They determine two torsion tensors $S$ and $ T $ defined by
 \begin{equation}\label{S and T}
    S(X,Y)=\tau(H\hat{X},H\hat{Y}) ,\qquad  T(\dot{X},Y)=\tau(V\hat{X},H\hat{Y}) ,
 \end{equation}
    where $ H\hat{X}\in H_zTM=ker\mu_z $ and $ V\hat{X}\in V_zTM=ker(\pi_{*})_z $, and $ (\pi_{*})_z $ is the tangent mapping of the canonical projection $ \pi $.
     There is a unique regular connection $ {}^c\nabla$ associated to the Finsler structure $F$ satisfying, $ {}^c\nabla_{\hat{Z}}g=0 $,  $ S(X,Y)=0 $ and $ g(\tau(V\hat{X},\hat{Y}),Z)=g(\tau(V\hat{X},\hat{Z}),Y) $, called the \textit{Cartan connection}. The condition ${}^c\nabla_{\hat{_Z}}g=0 $ is called \textit{metric compatibility} in both horizontal and vertical covariant derivatives, which is equivalent to
   \begin{equation}\label{compatibility}
   \hat{Z}g(X,Y)=g({}^c\nabla_{\hat{Z}}X,Y)+g(X,{}^c\nabla_{\hat{Z}}Y).
   \end{equation}
   It results from the last equation that the Cartan covariant derivative ${}^c\nabla$ is determined by the \textit{Finslerian Koszul formula}, see \cite{Ak3,BSY}.
  \begin{align}\label{kozol}
   2g({}^c\nabla_{\hat{X}}Y,Z)&=\hat{X}.g(Y,Z)+\hat{Y}.g(X,Z)-\hat{Z}.g(X,Y)+g(\tau(\hat{X},\hat{Y}),Z)
   +g(\tau(\hat{Z},\hat{X}),Y)\nn\\&+g(\tau(\hat{Z},\hat{Y}),X)+
    g(\varrho[\hat{X},\hat{Y}],Z)+g(\varrho[\hat{Z},\hat{X}],Y)+g(\varrho[\hat{Z},\hat{Y}],X).
    \end{align}
     According to the definition of connection $1$-form of Cartan connection we have $\omega^{i}_{j}:=\Gamma^{i}_{jk}dx^{k}+C^{i}_{jk}\delta y^{k},$
  where
 \begin{eqnarray*}
 \Gamma^{i}_{jk}:=\frac{1}{2}g^{il}(\delta _{j}g_{lk}+\delta_{k}g_{jl}-\delta _{l}g_{jk}),\quad  C^{i}_{jk}:=\frac{1}{2}g^{il}\dot{\partial}_{l}g_{jk},
\end{eqnarray*}
and
$ \delta _{i}:=\frac{\delta}{\delta x^{i}} $,
$ \dot{\p}_{i}:=\frac{\partial}{\partial y^{i}}.$
  Using $\hat{X}=H\hat{X}+V\hat{X}$, the Cartan covariant derivative is decomposed in the horizontal and vertical forms ${}^c\nabla_{\hat{X}}Y=\nabla_{H\hat{X}}Y+\nabla_{V\hat{X}}Y$.
In a local coordinate system the components of the Cartan connection  ${}^c\nabla_k$ are denoted here by
\begin{eqnarray*}%\label{Def;CartanCon}
{}^c\nabla_k=\nabla_{k}+ \dot{\nabla}_{k},
\end{eqnarray*}
wherein,$ \nabla_{k}:={}^c\nabla_{\frac{\delta}{\delta x^{k}} }$,
$ \dot{\nabla}_{k}:={}^c\nabla_{\frac{\partial}{\partial y^{k}}}$.
  Denote by $\Gamma^{i}_{jk}$ and $C^{i}_{jk}$  the \emph{horizontal} and the \emph{vertical}  coefficients of Cartan connection respectively. We have
\begin{eqnarray*}%\label{Def;CartanCon}
&\nabla_{k}\dot{\partial_{j}}=\Gamma^{i}_{jk}\dot{\partial _{j}},\qquad
& \dot{\nabla}_{k}\dot{\partial _{j}}=C^{i}_{jk}\dot{\partial_{j}},\nn \\
& \nabla_{k} \delta _{j} =\Gamma ^{i}_{jk}\delta_{i},\qquad
& \dot{\nabla}_{k} \delta _{j}=C^{i}_{jk}\delta _{i}.\nn
\end{eqnarray*}
In a local coordinate system,
the \textit{horizontal} and \textit{vertical}  Cartan covariant derivatives of an arbitrary $(1,1)$-tensor field on $\pi^*TM$ with the components $T_{i}^{j}$ are given by
\begin{eqnarray}\label{Def;hvCartanCon}\label{Def;CartanConn}
& \nabla_{k}T_{i}^{j}=\delta_{k}T_{i}^{j}-T_{r}^{j}\Gamma ^{r}_{ik}+T_{i}^{r}\Gamma ^{j}_{rk},\\
& \dot{\nabla}_{k}T_{i}^{j}=\dot{\partial }_{k}T_{i}^{j}-T_{r}^{j}C ^{r}_{ik}+T_{i}^{r}C ^{j}_{rk}.\nn
\end{eqnarray}
The components of Cartan hh-curvature tensor are given by
\begin{equation*} \label{77}
R^{i}_{jkm}=\delta_{k}\Gamma^{i}_{jm}-\delta_{m}\Gamma^{i}_{jk}+
\Gamma^{i}_{s k}\Gamma^{s}_{jm}-\Gamma^{i}_{s m}\Gamma^{s}_{jk}
+R^{s}_{km}C^{i}_{s j},
\end{equation*}
where, $R^{i}_{km}:=y^{p} R^{i}_{ pkm}$.
For a non-null $ y\in T_xM $,  the trace of  hh-curvature is called  \textit{Riemann curvature}. It is given by
$ R_y(u)=R^i_ku^k\frac{\partial}{\partial x^i},$     where
 \begin{equation*}
  R^i_k(y):=\frac{\partial G^i}{\partial x^k}-1/2\frac{\partial^2G^i}{\partial y^k\partial x^j}y^j+G^j\frac{\partial^2G^i}{\partial y^k\partial y^j}-1/2\frac{\partial G^i}{\partial y^j}\frac{\partial G^j}{\partial y^k}.
 \end{equation*}
 The \textit{Ricci scalar} is defined by $ Ric:=R^i_i $, see \cite[p. 331]{BCS}. Here, we use Akbar-Zadeh's definition of \textit{Ricci tensor} as follows
  $Ric_{ik}:=1/2(F^2Ric)_{y^iy^k} $, see \cite{Ak3}. Let $ N^i_j=\frac{\partial G^i}{\partial y^j},$ {}  $  l^i=\frac{y^i}{F}$ and $ \hat{l}=l^i\frac{\delta}{\delta x^i}=l^i(\frac{\partial}{\partial x^i}-G^k_i\frac{\partial}{\partial y^k}) $.  By homogeneity we have $ Ric_{ij} l^il^j=Ric $. Let $\tilde{F} $ be another Finsler structure on $M$.
  In this paper we deal with the forward geodesics and  the word ``geodesic"    refers to the forward geodesic.
  If any geodesic of $(M,F)$  coincides with a geodesic of $(M,\tilde{F})$ as set of points and vice versa, then the change  $F\rightarrow \tilde{F}$ of the metric is called \emph{projective } and $F$ is said to be \emph{projective} to $\tilde{F}$.
   A Finsler space $(M,F)$  is projective to another Finsler space $(M,\tilde{F})$, if and only if there exists a  1-homogeneous scalar field $p(x,y)$ satisfying
       $ \tilde{G}^i(x , y)=G^i(x , y)+p(x,y)y^i$.
      The scalar field $p(x,y)$ is called the \emph{projective factor} of the projective change under consideration.
    For a tangent plane $ P\subset T_pM $ and a non-zero vector $ y\in T_pM $, the \textit{flag curvature} $ K(P,y) $ is defined by \[K(P,y)=\frac{g_y(u,R_y(u))}{g_y(y,y)g_y(u,u)-g_y(y,u)^2},\] where $ P=span\{y,u\} $.
   When $F$ is Riemannian, $K(P,y)=K(P)$ is independent of $y\in P$ and is just the sectional curvature in Riemannian geometry.
   We say that $ F $ is of \textit{scalar curvature} if for any $ y\in T_pM $, the flag curvature $ K(P,y)=K(y) $ is independent of $ P $ containing $ y\in T_pM $.
  In a local coordinate system $(x^i,y^i)$ on $TM$, this is equivalent to saying
 \[R^i_k=KF^2\{\delta^i_k-F^{-1}F_{y^k}y^i\}.\]
 If $ K $ is constant, then $ F $ is said to be of \textit{constant flag curvature}.
 A Finsler metric $ F $ that satisfies the relation, $Ric=(n-1)K(x)$, for some functions $ K $ on $ M $, is called \textit{Einstein metric}. It is well known the relation between the Ricci scalar $ Ric $ and the Ricci tensor $ Ric_{ij} $, tells us,
 \begin{equation}\label{Ric}
 Ric=(n-1)K(x)  \Longleftrightarrow Ric_{ij}=(n-1)K(x)g_{ij}.
 \end{equation}
 See \cite{BR}. If the function $K$ is constant, then $F$ is called \textit{Ricci-constant}.
  A Finsler space is called a \textit{Berwald space } if the Berwald connection coefficients, namely $ (G^i)_{y^jy^k}$, do not depend on $ y$.
   In particular, all Riemannian and locally Minkowskian spaces are Berwaldian, see \cite{BR}.
The Finsler structure of a \textit{Randers metric} on the smooth n-dimensional manifold $ M $ is given by $ F=\alpha+\beta $, where $\alpha(x,y):=\sqrt{a_{ij}y^iy^j},$  is a Riemannian metric and $\beta(x,y):=b_i(x)y^i,$ is a 1-form, see \cite{ChSh}. 
We will need the following proposition in the sequel.
\begin{pro}\cite{BR}\label{p1}
Let $ (M,F) $ be a connected compact boundaryless Einstein Randers manifold with constant Ricci scalar $ Ric $.\\
\textbullet If $ Ric<0 $, then $ (M,F) $ is Riemannian.\\
\textbullet If $ Ric=0 $, then $ (M,F) $ is Berwaldian.
\end{pro}
\subsection{Geodesics and circles on Finsler manifolds}
We recall here a natural definition given in \cite{BS} of a circle in a Finsler manifold.
Let $c:I\subset\mathbb{R}\longrightarrow M$ be a smooth curve parameterized by the arc length $s$ on a Finsler manifold $(M,F)$.
 Consider a unitary normal vector field $Y$ along $c$ and a positive constant $\kappa$ such that $\nabla_{\dot{c}}X=\kappa Y$ and
$\nabla_{\dot{c}}Y=-\kappa X,$ where, $X:=\dot{c}=\frac{dc}{ds}$ is the unitary tangent vector field at each point $c(s)$ and $\nabla_{\dot{c}}$ is the Cartan covariant derivative along $c$.
 The numbers $\kappa$ and $\frac{1}{\kappa}$ are called curvature and radius of the circle, respectively. A \emph{geodesic circle} on a Finsler space $(M,F)$ is defined to be a smooth curve $c: I\longrightarrow M$ for which the first Frenet curvature $\kappa_1:=\kappa$, is constant and the second Frenet curvature $\kappa_2$, vanishes identically. That is, $\frac{d\kappa_1}{ds}=0$ and $\kappa_{2}=0$, see, \cite{SedBid,BS}. If in the definition of a geodesic circle we exclude geodesic or equivalently the trivial case, $\kappa_{1}=0$, then we obtain the definition of a circle on a Finsler space.
 A  conformal change of metric is said to be \textit{concircular} if it maps geodesic circles into geodesic circles. The following theorems will be used in the sequel.
    \begin{theo}\cite{Bi2}\label{pre}
    Let $ (M,F) $  be a Finsler manifold. A necessary and sufficient condition for a conformal change $ \bar{g}=e^{2\varphi(x)}g $ to be  concircular, is the function $ \varphi $ be a solution of the partial differential equation
    \bes%\label{a9}
      {}^c\nabla_i\varphi_j - \varphi_i \varphi_j=\Phi g_{ij},
   \ees
     where $ \varphi_j=\partial\varphi/\partial x^j $, ${}^c\nabla_i $ is the Cartan horizontal derivative and $ \Phi $ is a certain scalar function.
    \end{theo}
    \begin{theo}\cite{SedBid}\label{Th;p1}
If a forward geodesically complete Finsler manifold $ (M,F) $ admits a circle preserving change of metric, then the indicatrix is conformally diffeomorphic to the Euclidean sphere $ S^{n-1} $ in $ \mathbb{R}^n $.
\end{theo}
\begin{theo}\cite{SedBid}\label{Th;p2}
Let $ (M,F) $ be a forward geodesically complete absolutely homogeneous Finsler manifold of scalar flag curvature. If $ (M,F) $ admits a nontrivial circle preserving change of metric, then it is a Riemannian manifold of constant sectional curvature.
\end{theo}
    \section{Conformal change of Cartan connection}
\subsection{Conformal transformations on Finsler manifolds}
Let $F$ and $\bar{F}$ be two Finsler structures on an n-dimensional manifold $M$. A diffeomorphism $f:(M,F)\to(M,\bar{F})$ is called \textit{conformal transformation}  or  simply a \textit{conformal change} of metric, if and only if there exists a scalar function $\varphi(x)$ on $M$ such that $\bar{F}(x,y)=e^{\varphi(x)}F(x,y)$.
 One can easily show that, the function $\varphi(x,y)$ is independent of the direction $y$, or equivalently $\frac{\p \varphi}{\p y^i}=0.$
Here, the terms mapping and transformation will be used
interchangeably.
Assuming $\bar{F}(x,y)=e^{\varphi(x)}F(x,y)$ the above relation becomes\begin{equation}\label{z8}
\bar{g}=e^{2\varphi(x)}g,
\end{equation}
where, $\bar{g}:=f^*g$.
The diffeomorphism $f$ is said to be \textit{homothetic} if $\varphi$ is constant and \textit{isometric} if $ \varphi $ vanishes in every point of $M$.

Throughout this article, the objects of $({M}, \bar{F})$ will be shown with a bar and we shall always assume that the line elements $(x,y)$ and $(\bar{x},\bar{y})$ on $(M,F)$ and $({M},\bar{F})$ are chosen such that $ \bar{x^i}=x^i$ and $ \bar{y^i}=y^i $ holds, unless a contrary assumption is explicitly made. If we show the corresponding Finsler metric tensors by $ g $ and $ \bar{g} $, then Eq. \eqref{z8} is written in the following local forms:
\begin{equation}\label{z9}
\bar{g}_{ij}(x,y)=e^{2\varphi(x)}g_{ij}(x,y), \quad \bar{g^{ij}}(x,y)=e^{-2\varphi(x)}g^{ij}(x,y),
\end{equation}
where $ g^{ij} $ is the inverse matrix defined by $ g_{ij}g^{ik}=\delta^k_j $. Eq. \eqref{z8} and definition of \textit{Cartan tensor} yield
\begin{equation}\label{Cbar}
\bar{C}^i_{jk}(x,y)=C^i_{jk}(x,y),\quad \bar{C}_{ijk}=e^{2\varphi}C_{ijk},
\end{equation}
 where $ C^i_{jk}=g^{il}C_{ljk}=1/2g^{il}\frac{\partial g_{lj}}{\partial y^k}$.

   It is well known that after a conformal change of metric, the Christoffel symbols $\gamma^i_{jk}$, as a geometric object on a Finsler manifold $(M,F)$, satisfies
\begin{equation}\label{chris}
\bar{\gamma}^i_{jk}=\gamma^i_{jk}+(\delta^i_j\delta^h_k+\delta^i_k\delta^h_j-g^{ih}g_{jk})\varphi_h,
\end{equation}
where, $\varphi_h=\partial\varphi/\partial x^h$ and $\delta^i_j$ is the Kronecker delta, see for instance, \cite[page 28]{Ha}.
 Contracting both side of \eqref{chris} by $ y^jy^k $ and using $ 2G^i=\gamma^i_{jk}y^jy^k $, we obtain the relation between $ G^i $ and $ \bar{G^i} $ in conformal Finsler spaces as follows.
   \begin{equation}\label{conf}
   \bar{G}^i=G^i-B^{ir}\varphi_r,
   \end{equation} where $ B^{ir}:=(\frac{F^2}{2}g^{ir}-y^ry^i)  $.\\
    By differentiation of \eqref{conf} with respect to $ y^j $ we have
    \begin{equation}\label{Gijbar}
    \bar{G}^i_j=G^i_j-B^{ir}_j\varphi_r,
    \end{equation}
    where $ B^{ir}_j=y_jg^{ir}-F^2C^{ir}_j-\delta^r_jy^i-y^r\delta^i_j  $ and $ C^{ir}_j:=g^{is}C^r_{sj}=-\frac{1}{2}\frac{\partial g^{ir}}{\partial y^j} $, for more details see \cite{Ha}.\\
    By definition of $ \frac{\delta}{\delta x^k} $ and \eqref{Gijbar} we get
    \begin{equation*}
    \bar{\frac{\delta}{\delta x^k}}=\frac{\partial}{\partial x^k}-\bar{G}^i_k\frac{\partial}{\partial y^i}=\frac{\partial}{\partial x^k}-(G^i_k-B^{ir}_k\varphi_r)\frac{\partial}{\partial y^i}=\frac{\delta}{\delta x^k}+B^{ir}_k\varphi_r\frac{\partial}{\partial y^i}.
    \end{equation*}
    If we put $L_k:=B^{ir}_k\varphi_r\frac{\p}{\p y^i}$, then the above relation becomes
    \begin{equation}\label{deltabar}
    \bar{\frac{\delta}{\delta x^k}}=\frac{\delta}{\delta x^k}+L_k.
    \end{equation}
   %Analogous to the Riemannian geometry, by straightforward calculation we have the following results in Finsler geometry.
   % \cite{Yam}, \cite[p. 6] {Bi3}.
   The following identities are well known;
      \begin{itemize}
   \item[(1)]$ [\frac{\delta}{\delta x^i},\frac{\delta}{\delta x^j}]=R^h_{ij}\frac{\partial}{\partial y^h} $,
   \item[(2)]$ [\frac{\delta}{\delta x^i},\frac{\partial}{\partial y^j}]=\frac{\partial N^h_{i}}{\partial y^j}\frac{\partial}{\partial y^h}=N^h_{ij}\frac{\partial}{\partial y^h} $,
   \item[(3)] $ [\frac{\partial}{\partial y^i},\frac{\partial}{\partial y^j}]=0. $
   \end{itemize}
   Therefore by definition of the function $\varrho:TTM_0\to \pi^*TM$ and the above identities, we have
   \begin{equation}\label{roh.braket}
   \varrho([\frac{\delta}{\delta x^i},\frac{\delta}{\delta x^j}])=0,\quad \varrho([\frac{\delta}{\delta x^i},\frac{\partial}{\partial y^j}])=0,\quad \varrho([\frac{\partial}{\partial y^i},\frac{\partial}{\partial y^j}])=0.
   \end{equation}
    \bpr\label{lem;nab-bar}
    Let $ (M,F) $ and $ (M,\bar{F}) $ be two conformal Finsler manifolds with the Cartan connections $ {}^c\nabla $ and $ {}^c\overline{\nabla} $ respectively. The related covariant derivatives satisfy
\begin{align}\label{nab-bar}
{}^c\overline{{\nabla}}_{\hat{X}}Y=&{{}^c\nabla}_{\hat{X}}Y+(H\hat{X}\varphi)Y+(H\hat{Y}\varphi)X
-(\nabla\varphi )g(X,Y)+\notag\\&T(L_i,Y)+T(L_j,X)-g(T(L_t,Y),X)\frac{\partial}{\partial x^s}g^{ts}.
\end{align}
 \epr
\begin{proof}
  Let us assume $\hat{X}=\frac{\delta}{\delta x^i}$, $\hat{Y}=\frac{\delta}{\delta x^j}$ and $\hat{Z}=\frac{\delta}{\delta x^k}$,  by definition of $\varrho$ we have $\varrho(\hat{X})=X=\frac{\p}{\p x^i}$, $Y=\frac{\partial}{\partial x^j}$ and $Z=\frac{\partial}{\partial x^k}$. On the other hand by the hh-torsion freeness of Cartan connection we have $\tau(\frac{\delta}{\delta x^i},\frac{\delta}{\delta x^j})=0$, $\tau(\hat{X},\hat{Y})=\tau(\hat{Z},\hat{Y})=\tau(\hat{X},\hat{Z})=0$. Therefore from \eqref{roh.braket} the  Koszul formula \eqref{kozol} reduces to
   \begin{equation}\label{nabla-horizontal}
   2g(\nabla_{\frac{\delta}{\delta x^i}}\frac{\partial}{\partial x^j},\frac{\partial}{\partial x^k})=\frac{\delta}{\delta x^i}g_{jk}+\frac{\delta}{\delta x^j}g_{ik}-\frac{\delta}{\delta x^k}g_{ij}.
   \end{equation}
   Rewriting \eqref{nabla-horizontal} for $(M,\bar{g})$ and using \eqref{deltabar} and $L_i e^{2\varphi}=0$ we have
   \begin{align*}
   2\bar{g}(\bar{\nabla}_{\frac{\bar{\delta}}{\delta x^i}}\frac{\partial}{\partial x^j},\frac{\partial}{\partial x^k})=&\frac{\bar{\delta}}{\delta x^i}\bar{g}_{jk}+\frac{\bar{\delta}}{\delta x^j}\bar{g}_{ik}-\frac{\bar{\delta}}{\delta x^k}\bar{g}_{ij}\\=&(\frac{\delta}{\delta x^i}+L_i)(e^{2\varphi}g_{jk})+(\frac{\delta}{\delta x^j}+L_j)(e^{2\varphi} g_{ik})-(\frac{\delta}{\delta x^k}+L_k)(e^{2\varphi}g_{ij})\\
   =&e^{2\varphi}(\frac{\delta}{\delta x^i}g_{jk}+\frac{\delta}{\delta x^j} g_{ik}-\frac{\delta}{\delta x^k}g_{ij})+2e^{2\varphi}(\frac{\partial\varphi}{\partial x^i}g_{jk}+\frac{\partial\varphi}{\partial x^j}g_{ik}-\frac{\partial\varphi}{\partial x^k}g_{ij})\\
    &+L_i(e^{2\varphi}g_jk)+L_j(e^{2\varphi}g_{ik})-L_k(e^{2\varphi}g_{ij})\\
   =& 2e^{2\varphi}g(\nabla_{\frac{\delta}{\delta x^i}}\frac{\partial}{\partial x^j},\frac{\partial}{\partial x^k})+2e^{2\varphi}(\frac{\partial\varphi}{\partial x^i}g_{jk}+\frac{\partial\varphi}{\partial x^j}g_{ik}-\frac{\partial\varphi}{\partial x^k}g_{ij})\\&+ e^{2\varphi}(L_i(g_{jk})+L_j(g_{ik})-L_k(g_{ij})).
   \end{align*}
   Replacing $ \bar{g}=e^{2\varphi}g $ in the left hand side of the above relation we get
   \begin{align}\label{nablabarhor}
   {g}(\bar{\nabla}_{\frac{\bar{\delta}}{\delta x^i}}\frac{\partial}{\partial x^j},\frac{\partial}{\partial x^k})=&g(\nabla_{\frac{\delta}{\delta x^i}}\frac{\partial}{\partial x^j},\frac{\partial}{\partial x^k})+\varphi_ig_{jk}+\varphi_jg_{ik}-\varphi_kg_{ij}\notag\\&+\frac{1}{2}(L_i(g_{jk})
   +L_j(g_{ik})-L_k(g_{ij})).
   \end{align}
   On the other hand
   \begin{equation*}
   L_i(g_{jk})=B^{mr}_i\varphi_r\frac{\partial}{\partial y^m}g_{jk}=2B^{mr}_i\varphi_rC_{mjk}=2B^{mr}_i\varphi_rC^l_{mj}g_{kl}.
   \end{equation*}
     And
   \[L_k(g_{ij})=2B^{mr}_k\varphi_rC^l_{mj}g_{il}=2B^{mr}_t\varphi_rC^l_{mj}g_{il}\delta^t_k=
   2B^{mr}_t\varphi_rC_{mji}g^{ts}g_{sk}.\]
   Therefore \eqref{nablabarhor} becomes
   \begin{align*}
      {g}(\bar{\nabla}_{\frac{\bar{\delta}}{\delta x^i}}\frac{\partial}{\partial x^j},\frac{\partial}{\partial x^k})=&g(\nabla_{\frac{\delta}{\delta x^i}}\frac{\partial}{\partial x^j},\frac{\partial}{\partial x^k})+\varphi_ig_{jk}+\varphi_jg_{ik}-\varphi_kg_{ij}+\notag\\&B^{mr}_i\varphi_rC^l_{mj}g_{kl}
      +B^{mr}_j\varphi_rC^l_{mi}g_{kl}-B^{mr}_t\varphi_rC_{mji}g^{ts}g_{sk}.
 \end{align*}
   Hence
   \begin{align}\label{nablabarhor1}
   \bar{\nabla}_{\frac{\bar{\delta}}{\delta x^i}}\frac{\partial}{\partial x^j}=&{\nabla}_{\frac{{\delta}}{\delta x^i}}\frac{\partial}{\partial x^j}+\varphi_i\frac{\partial}{\partial x^j}+\varphi_j\frac{\partial}{\partial x^i}-(\nabla\varphi)g_{ij}+B^{mr}_i\varphi_rC^l_{mj}\frac{\partial}{\partial x^l}\notag\\
   &+B^{mr}_j\varphi_rC^l_{mi}\frac{\partial}{\partial x^l}-B^{mr}_t\varphi_rC_{mji}g^{ts}\frac{\partial}{\partial x^s},
   \end{align}
   where, $\nabla\varphi$ is the gradient of $ \varphi $ and we have $\varphi_k=g(\nabla\varphi,\frac{\partial}{\partial x^k})$.\\
   %It should be noted that the above relation is similar to the usual relation between $ \Gamma^k_{ij} $ and $ \bar{\Gamma}^k_{ij} $ found in \cite[p. 31]{Has} for Cartan connection.\\
   By the definition of  torsion   \eqref{S and T}, we have
   $ T(L_i,\frac{\partial}{\partial x^j})=\nabla_{L_i}\frac{\partial}{\partial x^j}=B^{mr}_i\varphi_rC^l_{mj}\frac{\partial}{\partial x^l} $, hence \eqref{nablabarhor1} becomes
   \begin{align}\label{nablahor2}
   \bar{\nabla}_{\frac{\bar{\delta}}{\delta x^i}}\frac{\partial}{\partial x^j}=&{\nabla}_{\frac{{\delta}}{\delta x^i}}\frac{\partial}{\partial x^j}+\varphi_i\frac{\partial}{\partial x^j}+\varphi_j\frac{\partial}{\partial x^i}-\nabla\varphi g_{ij}+\notag\\&T(L_i,\frac{\partial}{\partial x^j})+T(L_j,\frac{\p }{\p x^i})-g(T(L_t,\frac{\partial}{\partial x^j}),\frac{\partial}{\partial x^i})\frac{\partial}{\partial x^s}g^{ts}.
   \end{align}
  % where, $ (\nabla\varphi)^T $ is the tangential part of the gradient of $ \varphi $.\\
   Next assume that $\hat{X}=\frac{\partial}{\partial y^i}$ is the vertical and $\hat{Y}=\frac{\delta}{\delta x^j}$ and $\hat{Z}=\frac{\delta}{\delta x^k}$ are the  horizontal derivatives,  we have $ X=\varrho(\hat{X})=0$, $Y=\frac{\partial}{\partial x^j} $ and $ Z=\frac{\partial}{\partial x^k}$.  Replacing these values in the Koszul formula \eqref{kozol} we get
   \begin{equation}\label{z1}
   2g(\nabla_{\frac{\partial}{\partial y^i}}\frac{\partial}{\partial x^j},\frac{\partial}{\partial x^k})=\frac{\partial}{\partial y^i}g_{jk}+g(\tau(\frac{\partial}{\partial y^i},\frac{\delta}{\delta x^j}),\frac{\partial}{\partial x^k})+g(\tau(\frac{\delta}{\delta x^k},\frac{\partial}{\partial y^i}),\frac{\partial}{\partial x^j}).
   \end{equation}
   By the anti-symmetric property of torsion we have
   \[g(\tau(\frac{\partial}{\partial y^i},\frac{\delta}{\delta x^j}),\frac{\partial}{\partial x^k})=g(\tau(\frac{\partial}{\partial y^i},\frac{\delta}{\delta x^k}),\frac{\partial}{\partial x^j})=-g(\tau(\frac{\delta}{\delta x^k},\frac{\partial}{\partial y^i}),\frac{\partial}{\partial x^j}),\]
  hence
  \eqref{z1} yields
  \begin{equation}\label{nablavertic}
  2g(\nabla_{\frac{\partial}{\partial y^i}}\frac{\partial}{\partial x^j},\frac{\partial}{\partial x^k})=\frac{\partial}{\partial y^i}g_{jk}=2C_{ijk}.
  \end{equation}
  %$g(\frac{\partial}{\partial x^j},\frac{\partial}{\partial x^k})=g_{jk}$,
  Similarly for $(M,\bar{F})$,   $\bar{g}(\frac{\p}{\p x^j},\frac{\p}{\p x^k})=\bar{g}_{jk}$ and we have
  \begin{equation}\label{nablabarvertic1}
  2\bar{g}(\bar{\nabla}_{\frac{\partial}{\partial y^i}}\frac{\partial}{\partial x^j},\frac{\partial}{\partial x^k})=2\bar{g}(\bar{C}^l_{ij}\frac{\partial}{\partial x^l},\frac{\partial}{\partial x^k})=2\bar{C}^l_{ij}\bar{g}_{lk}=2\bar{C}_{ijk}.
  \end{equation}
  Replacing $\bar{g}=e^{2\varphi}g$ in the left side of  \eqref{nablabarvertic1} and using \eqref{Cbar} we get
 \begin{equation}\label{nablabarver11}
 {g}(\bar{\nabla}_{\frac{\partial}{\partial y^i}}\frac{\partial}{\partial x^j},\frac{\partial}{\partial x^k})=C_{ijk}.
 \end{equation}
    From \eqref{nablavertic} and \eqref{nablabarver11} we get
  \begin{equation*}
  g(\bar{\nabla}_{\frac{\partial}{\partial y^i}}\frac{\partial}{\partial x^j},\frac{\partial}{\partial x^k})=g({\nabla}_{\frac{\partial}{\partial y^i}}\frac{\partial}{\partial x^j},\frac{\partial}{\partial x^k}).
  \end{equation*}
  Therefore we have
  \begin{equation}\label{nablabarver2}
  \bar{\nabla}_{\frac{\partial}{\partial y^i}}\frac{\partial}{\partial x^j}={\nabla}_{\frac{\partial}{\partial y^i}}\frac{\partial}{\partial x^j}.
  \end{equation}
    The decomposition $\hat{X}=H\hat{X}+V\hat{X}$ yields ${}^c\overline{\nabla}_{\hat{X}}Y=\overline{\nabla}_{H\hat{X}}Y+\overline{\nabla}_{V\hat{X}}Y $.
     Using \eqref{nablahor2} and \eqref{nablabarver2} we have the proof.
       \end{proof}
   Recall that on a Riemannian manifold,  the metric $g$ is independent of the direction $y$ and $L_k=B^{ir}_k\varphi_r\frac{\partial}{\partial y^i}$, hence we have $L_i=0$.
   Therefore, \eqref{nab-bar} reduces to the following relation, see \cite{OS}.
       \begin{align*}
     \overline{\nabla}_XY=\nabla_XY+(X\varphi)Y+(Y\varphi)X-(\nabla\varphi)g(X,Y), \qquad\forall X,Y\in {\cal X}(M).
    \end{align*}
 \section{Hessian and Laplacian on Finsler manifolds}
 Here based on the global Cartan connection the Hessian and Laplacian are defined.
 The natural definitions of Hessian and Laplacian considered here are in some senses more general than those given in \cite{Ak1, BSY,Sh,WX}  and contains some of them in special cases.
\subsection{Horizontal and vertical Hessian}
 Let $f\in C^\infty(M)$,  the \textit{Hessian} of $f$ in the Cartan connection ${}^c\nabla_{\hat{X}}$ is defined by
 \beq\label{Hess}
 Hess&:&\Gamma(TM_0)\times \Gamma(\pi^*TM) \to C^\infty (TM_0) \nn \\
  Hess(f)(\hat{X},Y)&=&g(Y , {}^c\nabla_{\hat{X}}(\nabla f)),
 \eeq
% \beq\label{Hess}
% Hess &:&(T_zTM_0)\times {\cal X}(M)\to C^\infty (TM_0) \nn \\
%  Hess(f)(\hat{X},Y)&=&(\hat{X}Y)f-({}^c\nabla_{\hat{X}}Y)f,
% \eeq
 for all $\hat{X}\in \Gamma(TM_0)$ and $Y\in\Gamma(\pi^*TM)$. %see \cite[p. 181]{WX}.\\
 Using \eqref{compatibility} the metric compatibility of Cartan connection for a gradient vector field, we have
 \beq\label{compatibility+}
   \hat{X}g(Y,\nabla f)&=&g({}^c\nabla_{\hat{X}}Y,\nabla f)+g(Y,{}^c\nabla_{\hat{X}}\nabla f)\nn\\
  (\hat{X} Y) f &=& {}^c\nabla_{\hat{X}}Yf+g(Y,{}^c\nabla_{\hat{X}}\nabla f).\nn
   \eeq
   Replacing the last equation in \eqref{Hess}, we find a definition for the \emph{Hessian} of $f$ in Cartan connection as follows
 \beq\label{Hess+}
  Hess(f)(\hat{X},Y)=(\hat{X}Y)f-({}^c\nabla_{\hat{X}}Y)f.
 \eeq
Note that $Hess(f)$ can be split in \emph{horizontal} and \emph{vertical} parts as follows
\begin{align*}
Hess(f)(\hat{X},Y)=Hess^H_f(\hat{X},Y)+Hess^V_f(\hat{X},Y).
\end{align*}
In terms of the local frame fields $ \{\frac{\delta}{\delta x^j},\frac{\p}{\p y^j}\}$ and $ \{\frac{\p}{\p x^j}\} $ from \eqref{Hess+} we have
\bes%\label{hor-hess}
 Hess^H_f(\frac{\delta}{\delta x^i},\frac{\partial}{\partial x^j})=\frac{\delta}{\delta x^i}(\frac{\partial f}{\partial x^j})-(\nabla_{\frac{\delta}{\delta x^i}}\frac{\partial}{\partial x^j})f=\frac{\delta}{\delta x^i}(\frac{\p f}{\p x^j})-\Gamma^k_{ij}\frac{\partial f}{\partial x^k},
 \ees
  see \cite{BSY}.
Since $f$ is a function of $x$ alone, the above equation reduces to
\begin{equation}\label{loc.hor-hess}
 Hess^H_f(\frac{\delta}{\delta x^i},\frac{\partial}{\partial x^j})=\frac{\partial^2 f}{\partial x^i\partial x^j}-\Gamma^k_{ij}\frac{\partial f}{\partial x^k}.
\end{equation}
Following \eqref{Hess+},  the vertical part of $Hess(f)$ is written
\begin{equation*}%\label{ver-hess}
 Hess^V_f(\frac{\partial}{\partial y^i },\frac{\partial}{\partial x^j})=\frac{\p}{\p y^i}(\frac{\p f}{\p x^j})-(\nabla_{_\frac{\p}{\p y^i}}\frac{\p}{\p x^j})f=\frac{\p}{\p y^i}(\frac{\p f}{\p x^j})-C^k_{ij}\frac{\p f}{\p x^k}.
\end{equation*}
Again since $f$ is a function of $x$ alone, the vertical Hessian reduces to
\begin{equation}\label{loc.ver-hess2}
 Hess^V_f(\frac{\partial}{\partial y^i },\frac{\partial}{\partial x^j})=-C^k_{ij}\frac{\p f}{\p x^k}.
\end{equation}
Therefore for every function $f$ on $M$ by \eqref{loc.hor-hess} and \eqref{loc.ver-hess2} we have
\be\label{Def;Hess}
Hess(f)(\hat{X},Y)=\frac{\partial^2 f}{\partial x^i\partial x^j}-(\Gamma^k_{ij}+C^k_{ij})\frac{\partial f}{\partial x^k}.
\ee
\br
Recall that the Hessian on a Riemannian manifold is defined  by\\ $Hess(f)(X,Y)=(XY)f-(\nabla_{X}Y)f,$ for all $X,Y\in {\cal X}(M).$
In a local coordinate on a Riemannian manifold
\bes
Hess(f)({X},Y)=\frac{\p^2 f}{\p x^i\p x^j}-\gamma^k_{ij}\frac{\p f}{\p x^k},
\ees
 where $\gamma^k_{ji}$ are formal Christoffel symbols given by \eqref{Eq;ChristRie}.
 \er
 \subsection{Horizontal and vertical Laplacian}
   Generally, a  Laplacian is the trace of a Hessian, and in our setting the \textit{ horizontal Laplacian} $ \Delta^H f $ and the  \textit{vertical Laplacian} $\Delta^V f$ of a real smooth  function $f$ on $M$  are defined respectively by
     \begin{align*}
\Delta^H f&=trace_g (Hess^H_f(\hat{X},Y)),\\
\Delta^V f&=trace_g (Hess^V_f(\hat{X},Y)),
\end{align*}
 where $Y\in \Gamma(\pi^*TM),$ and $\hat{X}\in\Gamma(TM_0)$. Equivalently, in a local coordinate system, using \eqref{loc.hor-hess} and \eqref{loc.ver-hess2} the  horizontal and the vertical Laplacian of $f$ are given by
  \[ \Delta^H f=g^{ij}(\frac{\p^2 f}{\p x^i\p x^j}-\Gamma^k_{ij}\frac{\p f}{\p x^k}),\]
  where $f\in C^\infty(M)$ and
   $$\Delta^Vf=g^{ij}(-C^k_{ij}\frac{\partial f}{\partial x^k}),$$
   respectively.
   The above horizontal Laplacian is also used in \cite[p. 362]{Ak1} and \cite{BSY}.
       Recall that a Finsler metric $g_{ij}$ defines an inner product on the sections of $\pi^*TM$, hence one needs to consider the complete lift of a vector field on $M$ to introduce the concept of conformal vector fields on Finsler geometry. More intuitively let $ V=v^i\frac{\partial}{\partial x^i} $ be a vector field on the smooth manifold $M$.
     The \emph{complete lift} and the \emph{horizontal lift} of $V$ are two globally defined vector fields on $TM_0$ given by $\tilde{V}=v^i(x)\frac{\partial}{\partial x^i} +y^j(\frac{\partial v^i}{\partial x^j})\frac{\partial}{\partial y^i}$,
    and ${}^h\!\hat{V}=v^i(x)\frac{\delta}{\delta x^i}$, respectively.
\br
      Let $\varphi $ be a function of $x$ alone and  $X=X^i\frac{\partial}{\partial x^i}$, a vector field on the smooth manifold $M$, clearly we have
       \begin{equation*}%\label{X(phi)=hatX(phi)}
     X\varphi=\tilde{X}\varphi \ \textrm{and} \quad X\varphi={}^h\!\hat{X}\varphi.
        \end{equation*}
        where $\tilde{X}$ and ${}^h\!\hat{X}$ are the  complete lift and the horizontal lift of $X$, respectively.\\
         If necessary, we may restrict without loss of generality the vector field $\hat{X}$  to the complete lift $\tilde{X}$ on $TM_0$.
 \er
\section{Schwarzian of the conformal diffeomorphisms}
   The classical notions of Schwarzian derivative and Schwarzian operator of an analytic function on a plane, are generalized for a conformal mapping on a Riemannian manifold, by Osgood and Stowe, see \cite{OS}.
Here, inspiring the Riemannian Schwarzian operator developed in \cite{Ca, OS} and using the above definitions of the gradient, Hessian and Laplacian for Cartan connection on a Finsler manifold,  a natural definition of Schwarzian is given for conformal mappings.
   %%%%%%%%%%%%%%%%%%%%%%%%%%%
\bd\label{Def;SchTen}
Let $(M,F)$ and $(M,\bar{F})$ be two conformally related Finsler manifolds, where $\bar{F}=e^{\varphi}F$.
 The \emph{Schwarzian tensor} $B_{_F}(\varphi)$ is a  symmetric traceless $(0,2)$-tensor field defined by
\begin{align}\label{c2}
B_{_F}(\varphi)(\hat{X},Y)=\textrm{Hess}(\varphi)(\hat{X},Y)-(d\varphi\otimes d\varphi)(\varrho\hat{X},Y)-\frac{1}{n}(\Delta\varphi-\|grad\varphi\|^2)g(\varrho\hat{X},Y),
\end{align}
for all $\hat{X}\in \Gamma(TM_0)$ and $Y\in \Gamma(\pi^*TM)$  where   $\|grad\varphi\|^2=\varphi^i\varphi_i$, \, $\varphi^i=g^{ij}\varphi_j$  and
$g$ is the inner product on $\pi^*TM$ derived from the Finsler structure $F$.
\ed
  We define  the Schwarzian derivative of a conformal transformation as an operator applying to the vector fields $\hat X$ on $TM_0$ in the following sense.
   \bd\label{Def;SchDer}
 Let $F$ and $\bar{F}$ be two conformally related Finsler metrics on $M$.
      The \emph{Schwarzian derivative} of a conformal map $f:(M,F)\to (M,\bar{F})$ with $\bar{F}=e^{ \varphi}F$, at a point $ x\in M $,
       is a linear map
\begin{align*}
 S_F(f):\Gamma(TM_0)  &\longrightarrow \Gamma(\pi^*TM) ,\\
  S_F(f) \hat{ X}&= {}^c\nabla_{\hat{X}}(\nabla \varphi)-g(\nabla\varphi,\varrho\hat{X})\nabla\varphi-\frac{1}{n}(\Delta\varphi-\|grad\varphi\|^2) \varrho \hat{X},
 \end{align*}
    where  $ \hat{X}\in \Gamma(TM_0)$ and $\varrho \hat{X}=X$.
     \ed
     Lack of ambiguity, we denote simply the Schwarzian derivative $S_F(f)$  by $S(f)$.
 From which we obtain
   \begin{align}\label{Sf,P(phi)}
  g(S_{_F}&(f)(\hat{X}),Y)\nn\\
  &=g({}^c\nabla_{\hat{X}}\nabla \varphi,Y)-g(\nabla\varphi,\varrho\hat{X})g(\nabla\varphi,Y)
  -1/n(\Delta\varphi-\|\textrm{grad}\varphi\|^2)g(\varrho\hat{X},Y)\nn\\
 &=Hess(\varphi)(\hat{X},Y)-g(\nabla\varphi,{X})g(\nabla\varphi,Y)
 -1/n(\Delta\varphi-\|\textrm{grad}\varphi\|^2)g({X},Y)\nn\\
  &=Hess(\varphi)(\hat{X},Y)-({X}\varphi)(Y\varphi)
 -1/n(\Delta\varphi-\|\textrm{grad}\varphi\|^2)g({X},Y)\nn\\
 &=B_{_F}(\varphi)(\hat{X},Y),
    \end{align}
    where, $\hat{X}\in \Gamma(TM_0),$ $\varrho(\hat{X})=X$ and  $X,Y\in \Gamma(\pi^*TM)$.
     The the Schwarzian tensor $B_{_F}(\varphi)$ plays a similar role as the Schwarzian derivative of the conformal diffeomorphism $f$, in the above sense. The Schwarzian derivative  $S_F(f)$ in Definition\ref{Def;SchDer} applies on the vector fields  $\hat{X}\in \Gamma(TM_0)$, as well the Schwarzian tensor $B_{_F}(\varphi))(\ .\ ,Y)$ in  Definition\ref{Def;SchTen}.

       As a  raison d'\^{e}tre for our terminology, we determine $S(f)$, when %$f$ is an analytic function in the plane,
       $f$ is a real function on $\R$ with the Euclidean metric
       $ g=|dx|^2 $ or $ g_{ij}=\delta_{ij}$. In this case $\varphi=log|f'|$ and $ f^*|dx|=|f'||dx| $, by computing in standard coordinates one
          gets
        \be\label{1*}
        \quad e^\varphi=f'.
         \ee
         By differentiating  \eqref{1*} with respect to $x$, we have
         \be\label{2**}
        \varphi' e^\varphi=f''.
         \ee
        From \eqref{1*} and \eqref{2**} we get $\varphi'=\frac{f''}{f'}.$
        A second differentiation of \eqref{2**} leads $e^\varphi(\varphi''+\varphi'^2)=f''', $ and $ \varphi''+\varphi'^2=\frac{f'''}{f'}.\\
      $ Replacing the last equation in the definition of $ S(f) $ and using \eqref{1*} and \eqref{2**}, we have
      \begin{eqnarray*}
      S(f)&=&\frac{f'''}{f'}-\frac{3}{2}(\frac{f''}{f'})^2=(\varphi''+\varphi'^2)
      -\frac{3}{2}(\varphi'^2)\\
 &=& \varphi''- \varphi'^2+\frac{1}{2}   |\varphi'|^2\\
      %=\varphi''-\varphi'^2-(1/2    |\varphi'|^2-|\varphi'|^2),
 &=&Hess\varphi-d\varphi\otimes d\varphi -(\Delta\varphi-\|grad\varphi\|^2)g\\
 &=&B_F(\varphi),
       \end{eqnarray*}
        where for the dimension $n=1$, we have $Hess\varphi=\varphi'', \varphi'= d\varphi$ and $2\Delta\varphi= |\varphi'|^2$.
        As another verification of Definition \ref{Def;SchDer},  the replacement of the  Riemannian Hessian and Laplacian in this definition gives the well-known definition of the Riemannian Schwarzian derivative of a conformal map.
      In terms of a local coordinate system, using the local Hessian \eqref{Def;Hess}, we have
  \begin{align*}%\label{a8+}
       \big(B_F(\varphi)\big)_{ij}= \frac{\partial^2\varphi}{\partial x^i\partial x^j}-(\Gamma^h_{ij}+C^h_{ij})\varphi_h-\varphi_i\varphi_j-\frac{1}{n}
        (\Delta\varphi-\|grad\varphi\|^2)g_{ij},
  \end{align*}
        where, $\Gamma^h_{ij}$ and $C^h_{ij}$ are the Christoffel symbols of Cartan connection and the Cartan tensor, respectively and $\varphi_j=\frac{\partial\varphi}{\partial x^j}$.
         Using the Cartan horizontal and vertical covariant derivative formulas \eqref{Def;CartanConn}, and the fact that $\varphi$ is a function of $x$ alone, the last equation is written in the following familiar  form
  \begin{align}\label{a8}
        \big(B_F(\varphi)\big)_{ij}={}^c\nabla_i\varphi_j-\varphi_i\varphi_j
        -\frac{1}{n}(\Delta\varphi-\|grad\varphi\|^2)g_{ij},
  \end{align}
        where, ${}^c\nabla_i\varphi_j=:\frac{\partial^2\varphi}{\partial x^i\partial x^j} - ( \Gamma^h_{ij}+C^h_{ij})\varphi_h$ are the components of Cartan covariant derivative.

 We will refer to $ B_{_F}(\varphi) $ operator in the sequel, as the \emph{Schwarzian derivative}
        $ S_{_F}(f) $ of a conformal map $f:(M,F)\to (M,\bar{F})$ with $\bar{F}=e^{ \varphi}F$, in Finsler geometry.
   \section{Mobius mapping and Finsler manifolds}
   In  Riemannian geometry, a  conformal mapping is said to be Mobius if its Schwarzian derivative is zero.
   Therefore, we could consider a natural definition for  Mobius mappings on Finsler manifolds as follows.
 \bd
   A conformal diffeomorphism $f:(M,F)\to (M,\bar{F}),$ is called a \emph{Mobius mapping}, if the Schwarzian derivative $S_{_F}(f)$ vanishes. %$ S_{_F}(f)=B_F(\varphi)=0 $.
 \ed
 By means of \eqref{Sf,P(phi)}, we know  $S_{_F}(f)=0$, if and only if  $B_F(\varphi)=0$, $\forall Y \in\pi^*TM$.

 Clearly every isometry or every  identity map $i_d:(M,F)\to (M,F)$ is a Mobius mapping with $\varphi(x)=0$.

  If we put $ \Phi=\frac{1}{n}(\Delta\varphi-\|grad\varphi\|^2) $, then \eqref{a8} becomes $ B_F(\varphi)={}^c\nabla_i\varphi_j-\varphi_i\varphi_j-\Phi g_{ij} $. Hence vanishing of the Schwarzian tensor $ B_F(\varphi)=0 $, is equivalent to
   \be\label{a9+}
      {}^c\nabla_i\varphi_j - \varphi_i \varphi_j=\Phi g_{ij}.
   \ee
  \br\label{Re;Mob=cir}
   The equation \eqref{a9+} and Theorem\ref{pre} characterize the Mobius mappings in the sense that, a conformal diffeomorphism is Mobius if and only if it preserves geodesic circles.
    \er
  A property of Mobius mappings is given in the following Theorem.
 \bth\label{th;Big(sum)}
   Let $(M,F)$ be a Finsler manifold and $\varphi$ and  $\sigma$ the two real smooth functions on $M$ such that  $\bar{F}=e^{\varphi}F$. We have
   \begin{equation*}%\label{th;Big(sum)}
    B_F(\varphi+\sigma)=B_F(\varphi)+B_{\bar{F}}(\sigma)+A(\sigma),
   \end{equation*}
   where, $A(\sigma)=T(L_i,Y)\sigma+T(L_j,X)\sigma-g(T(L_t,Y),X)\frac{\p}{\p x^s}g^{st}\sigma$ and $L_i:=B^{kr}_i\varphi_r\frac{\p}{\p y^k}$.
   \eth
   \bpf
   By means of the definition of Hessian with respect to the Finsler metric $\bar{g}$,  and  the fact
   $(H \hat{X}\varphi)=\frac{\delta}{\delta x^i}\varphi=X\varphi$, using \eqref{nab-bar} we have
      \begin{align}\label{hessian-bar}
   Hess_{\bar{g}}(\sigma)(\hat{X},Y)&=\hat{X}Y(\sigma)-\bar{\nabla}_{\hat{X}}Y(\sigma)\notag\\
   &=\hat{X}Y(\sigma)-{\nabla}_{\hat{X}}Y(\sigma)-X(\varphi)Y(\sigma)-Y(\varphi){X}(\sigma)+g(X,Y)\nabla \varphi(\sigma)\notag\\
   &\quad -T(L_i,Y)\sigma-T(L_j,X)\sigma+g(T(L_t,Y),X)\frac{\p}{\p x^s}g^{st}\sigma\notag\\
   &=Hess_g(\sigma)(\hat{X},Y)-{X}(\varphi)Y(\sigma)-Y(\varphi){X}(\sigma)
   +g(X,Y)g(\nabla\varphi,\nabla\sigma)\notag\\
   &\quad -T(L_i,Y) \sigma-T(L_j,X) \sigma+g(T(L_t,Y),X)\frac{\p}{\p x^s}g^{st} \sigma,
   \end{align}
  where,  $X,Y\in\Gamma(\pi^*TM)$ and $\hat{X}\in\Gamma(TM_0)$.
   Let us denote the gradient and the Laplacian with respect to $\bar{g}$ by $ \bar{\nabla} $ and $ \bar{\Delta}$ respectively. If  $\{e_1,e_2,...,e_n\}$ is a local orthonormal frame on $M$ with respect to $\bar{g}$, that is, $\bar{g}(e_i,e_j)=\delta_{ij}$ then $\{e^{-\varphi}e_1,e^{-\varphi}e_2,...,e^{-\varphi}e_n\}$ is a local orthonormal frame  with respect to $g$ on $M$. By definition $\forall i,j,\  T(L_i,e_j)={}^c\nabla_{L_i}e_j=0$, hence from \eqref{hessian-bar} we get
   \begin{align}\label{laplas-bar}
   \bar{\Delta}\sigma&=\sum_{i=1}^{n}Hess_{\bar{g}}(\sigma)(\hat{e_i},e_i)\notag\\
   &=\sum_{i=1}^{n}\{Hess_g(\sigma)(e^{-\varphi}\hat{e_i},e^{-\varphi}e_i)-2{e_i}(\varphi)e_i(\sigma)
   +g(e^{-\varphi}e_i,e^{-\varphi}e_i)g(\nabla\varphi,\nabla\sigma)\}\notag\\
  &=\sum_{i=1}^{n}Hess_g(\sigma)(e^{-\varphi}\hat{e_i},e^{-\varphi}e_i)
   -2\bar{g}(\bar{\nabla}\varphi,\bar{\nabla}\sigma)+ne^{-2\varphi}{g}(\nabla\varphi,\nabla\sigma)\notag\\
 &=e^{-2\varphi}\Delta\sigma-2\bar{g}(\bar{\nabla}\varphi,\bar{\nabla}\sigma)
   +ne^{-2\varphi}{g}(\nabla\varphi,\nabla\sigma),
\end{align}
    where $ e_i\in\Gamma(\pi^*TM) $ and $ \hat{e_i} $ is its complete lift to $ \Gamma(TM_0) $.
   By means of \eqref{z9} and the definition of gradient of the scalar function $\varphi$ with respect to $\bar{g}$ we have
   \[\bar{\nabla}\varphi=\nabla_{\bar{g}}\varphi=\bar{g}^{ij}\frac{\partial\varphi}{\partial x^i}\frac{\partial}{\partial x^j}=e^{-2\varphi}g^{ij}\frac{\partial\varphi}{\partial x^i}\frac{\partial}{\partial x^j}=e^{-2\varphi}\nabla_g\varphi.\]
  Similarly $ \bar{\nabla}\sigma=\nabla_{\bar{g}}\sigma=e^{-2\varphi}\nabla_g\sigma $. Hence \eqref{laplas-bar} becomes
   \begin{align}\label{laplac-bar2}
   \bar{\Delta}\sigma=e^{-2\varphi}\{\Delta_g\sigma+(n-2)g(\nabla\varphi,\nabla\sigma)\}.
   \end{align}
   On the other hand
     \[\|\bar{\nabla}\varphi\|^2_{\bar{g}}=\bar{g}^{ij}\varphi_i\varphi_j
     =e^{-2\varphi}g^{ij}\varphi_i\varphi_j=e^{-2\varphi}\|\nabla\varphi\|^2_g,\]
   and  similarly $\|\bar{\nabla}\sigma\|^2_{\bar{g}}=e^{-2\varphi}\|\nabla\sigma\|^2_g  $.
   The equations \eqref{hessian-bar} and \eqref{laplac-bar2} and definition of Schwarzian tensor \eqref{c2} with $\varrho \hat{X}=X$ imply
     \begin{align*}
   B_{\bar{F}}(\sigma)(\hat{X},Y)&=Hess_{\bar{g}}(\sigma)(\hat{X},Y)
   -{X}(\sigma)Y(\sigma)-\frac{1}{n}\{\bar{\Delta}\sigma-\|\bar{\nabla}\sigma\|^2
 _{\bar{g}}\}\bar{g}({X},Y)\\
  &=Hess_g(\sigma)(\hat{X},Y)-{X}(\varphi)Y(\sigma)
   -Y(\varphi){X}(\sigma)+g(X,Y)g(\nabla\varphi,\nabla\sigma)\\
&\ \ -{X}(\sigma)Y(\sigma)-\frac{e^{-2\varphi}}{n}\{\Delta_g\sigma+(n-2)g(\nabla\varphi,\nabla\sigma)
-\|\nabla\sigma\|^2_g\}e^{2\varphi}g({X},Y)\\
&\ \ -T(L_i ,Y)\sigma-T(L_j,X)\sigma+g(T(L_t,Y),X)\frac{\p}{\p x^s}g^{st}\sigma\\
&=B_{F}(\sigma)(\hat{X},Y)-{X}(\varphi)Y(\sigma)-Y(\varphi){X}(\sigma)
   +\frac{2}{n}g(\nabla\varphi,\nabla\sigma)g({X},Y)\\
 &\ \ -T(L_i ,Y)\sigma-T(L_j,X)\sigma+g(T(L_t,Y),X)\frac{\p}{\p x^s}g^{st}\sigma.
 \end{align*}
If  we put $A(\sigma):= T(L_i ,Y)\sigma+T(L_j ,X)\sigma-g(T(L_t , Y),X)\frac{\partial}{\partial x^s}g^{st}\sigma$, the above equation yields
\begin{align}\label{Bg-bar}
B_{\bar{F}}(\sigma)=B_F(\sigma)-d\varphi\otimes d\sigma-d\sigma\otimes d\varphi+\frac{2}{n} g(\nabla\varphi ,\nabla\sigma)g-A(\sigma).
\end{align}
On the other hand
  \begin{align}\label{hessian-sum}
  Hess_g(\varphi+\sigma)-d(\varphi+\sigma)\otimes d(\varphi+\sigma)=&Hess_g(\varphi)+Hess_g(\sigma)-d\varphi\otimes d\varphi\notag\\
 &-d\varphi\otimes d\sigma-d\sigma\otimes d\varphi-d\sigma\otimes d\sigma.
  \end{align}
   Using the definition of Schwarzian tensor,  \eqref{1*} and \eqref{hessian-sum} we have
  \begin{align*}
  B_F(\varphi+\sigma)&= Hess_g(\varphi+\sigma)-d(\varphi+\sigma)\otimes d(\varphi+\sigma)-\frac{1}{n}\{\Delta(\varphi+\sigma)-\|\nabla(\varphi+\sigma)\|^2\}g\\&
  =Hess_g(\varphi)-d\varphi\otimes d\varphi-\frac{1}{n}\{\Delta\varphi-\|\nabla\varphi\|^2\}g
  \\ & \ \  +Hess_g(\sigma)-d\sigma\otimes d\sigma-\frac{1}{n}\{\Delta\sigma-\|\nabla\sigma\|^2\}g\\
  & \ \ -d\varphi\otimes d\sigma-d\sigma\otimes d\varphi+\frac{2}{n}g(\nabla\varphi,\nabla\sigma)g\\
  &=B_F(\varphi)+B_F(\sigma)-d\varphi\otimes d\sigma-d\sigma\otimes d\varphi+\frac{2}{n}g(\nabla\varphi,\nabla\sigma)g.
  \end{align*}
  The equation \eqref{Bg-bar} implies $B_{_F}(\varphi+\sigma)=B_{_F}(\varphi)+B_{_{\bar{F}}}(\sigma)+A(\sigma)$, and we have the proof.
   \epf
   \br
   If the conformal transformation is an identity or isometry, then we have $\varphi=0$ and Theorem\ref{th;Big(sum)}, yields
   \begin{equation}\label{Bg Bgbar}
    B_{_F}(\sigma)=B_{_{\bar{F}}}(\sigma)+A(\sigma).
    \end{equation}
   \er
   \begin{lem}\label{lem;Sch.of Compsition}
   Let $h:(M,F)\to(M,F')$ and $f:(M,F')\to(M,F'')$ be two conformal transformations on Finsler manifolds such that
   \be\label{eq;2conf}
   h^*g'=e^{2\varphi}g,\quad \textrm{and} \quad f^*g''=e^{2\sigma}g'.
   \ee
    Then
   \begin{equation}\label{composi}
    S_{_F}(f\circ h)=S_{_F}(h)+h^*S_{_F}(f).
    \end{equation}
\begin{proof}
Here, $\varphi$ and $\sigma$ are both real smooth functions on $M$.
By means of \eqref{eq;2conf} we have
\begin{align}\label{I}
(f\circ h)^*g''&=(h^*\circ f^*)g''=h^*(e^{2\sigma}g')=h^*(e^{2\sigma})h^*(g')=h^*(e^{2\sigma})e^{2\varphi}g\notag\\
&=(e^{2\sigma}\circ h)e^{2\varphi}g=e^{2(\sigma\circ h)}e^{2\varphi}g=e^{2(\varphi+(\sigma\circ h))}g.
\end{align}
By definition of Schwarzian derivative, if $f^*g''=e^{2\sigma}g' $, then $S_{F'}(f)=B_{F'}(\sigma)$, therefore from \eqref{I} we have $S_F(f\circ h)=B_F(\varphi+(\sigma\circ h))$. By Theorem\ref{th;Big(sum)} and the equations \eqref{eq;2conf} we get
\begin{align*}
S_F(f\circ h)=B_F(\varphi+(\sigma\circ  h))=B_F(\varphi)+B_{\bar{F}}(\sigma\circ h)+A(\sigma\circ h).
%&=B_{_F}(\varphi)+h^*B_{\bar{F}}(\sigma)+A(\sigma\circ h).
\end{align*}
From \eqref{Bg Bgbar} we have $B_F(\sigma \circ h)=B_{\bar{F}}(\sigma \circ h)+A(\sigma \circ h)$.
Therefore the above relation yields
$S_F(f \circ h)=B_F(\varphi)+B_F(\sigma \circ h) =B_F(\varphi)+h^*B_F(\sigma)=S_F(h)+h^*S_F(f).$
Therefore if $f$ and $h$ are both Mobius functions then their composition is also Mobius  and we have the proof.
\end{proof}
    \end{lem}
\emph{\textbf{Proof of Theorem \ref{th;group}.}}~
          Let $(M,F)$ and $(M,\bar{F})$ be two Finsler manifolds and  $f:(M,F)\to (M,\bar{F})$ a conformal transformation.
     In order to show the Mobius transformations from a group, we first show that the composites and inverses of Mobius transformations are Mobius.
      By means of Lemma \ref{lem;Sch.of Compsition} we have
      \[ 0= S_{_F}(i_d)=S_{_F}(f\circ f^{-1})=S_{_F}(f^{-1})+(f^{-1})^*S_{_F}(f), \]
     where $i_d$ is the identity map. Therefore, we have
$S_{_F}(f^{-1})=-(f^{-1})^*S_{_F}(f)$. Therefore, if $ f $ is a Mobius function, its inverse is  Mobius as well.
Next we show $S_{_F}(f)=S_{\bar{F}}(f).$ In fact, let $ \bar{F}=e^{\varphi}F$ and consider the composition  $f=h^{-1}\circ f\circ h$ given by $(M,{F})\xrightarrow{h}(M,{F})\xrightarrow{f}(M,\bar{F})\xrightarrow{h^{-1}}(M,\bar{F})$, where $h$  is the identity map and we have $S_{_F}(h)=S_{\bar{F}}(h^{-1})=0,$  hence by means of \eqref{composi} \[ S_{_F}(f)=S_{_F}((h^{-1}\circ f)\circ h)=S_F(h)+h^*S_{\bar{F}}(h^{-1}\circ f)=0+h^*\{S_{\bar{F}}(f)+f^*S_{\bar{F}}(h^{-1})\}.\]
   Since $ h $ is the identity map
   \[h^*(S_{\bar{F}}(f)+0)=S_{\bar{F}}(f).\]
   Hence  $ S_{_F}(f)=S_{\bar{F}}(f). $
Therefore, the set of Mobius transformations of $(M,F)$  forms a group.\\
 Let $ f:(M,F)\rightarrow(M,\bar{F}) $ be a homothety, since $ \varphi $ is constant, by definition  we get $ B_{_F}(\varphi)=0=S_{_F}(f).$ Thus any homothety, is a Mobius transformation and hence a subgroup of conformal transformations of $(M,F)$.
\hspace {\stretch{1}}$\Box$\\
\emph{\textbf{Proof of Theorem \ref{th;Con=Mobgroup}.}}
By means of \eqref{a9+}, Theorem\ref{th;group} and Theorem\ref{pre} we can easily see that a conformal diffeomorphism between two $n$-dimensional Finsler manifolds $(M,F)$ and $(M,\bar{F})$ is a Mobius transformation if and only if it maps all geodesic circles to geodesic circles. This completes the proof. For more details one can refer to \cite{BS}.
\hspace {\stretch{1}}$\Box$\\
Using the above properties of Schwarzian derivative we obtain proof of  the following rigidity theorems of complete Finsler manifolds.\\
\emph{\textbf{Proof of Theorems \ref{Th;MobiusComplet} and \ref{Th;MobiusScalar}.}}
By Remark\ref{Re;Mob=cir} the circle preserving or concircular maps are Mobius functions and vice versa. Therefore proof of these Theorems are a direct applications of Theorem\ref{Th;p1} and Theorem\ref{Th;p2}.
\hspace {\stretch{1}}$\Box$
\subsection{Schwarzian and Einstein Randers spaces}
As an application of the Schwarzian derivative, we can prove a rigidity theorem on Einstein Randers' spaces.\\
\emph{\textbf{Proof of Theorem \ref{Th;Randers}.}}
In \cite{BiS}, it is shown that in a Finsler manifold the projective parameter $p$ is a solution of the following ODE.
\begin{equation} \label{sch,Ric}
S(p(s))=\frac{\frac{d^3p}{ds^3}}{\frac{dp}{ds}}-\frac{3}{2}\Big[\frac{\frac{d^2p}{ds^2}}{\frac{dp}{ds}}\Big]^2 =\frac{2}{n-1}Ric_{jk}\frac{d{x}^j}{ds}\frac{d{x}^k}{ds}=\frac{2}{n-1}F^2Ric,
\end{equation}
where $ S(p) $ is the Schwarzian of \textquotedblleft $ p $\textquotedblright and \textquotedblleft $ s $\textquotedblright is the arc length parameter of a geodesic $ \gamma $.
 The projective parameter is unique up to a linear fractional transformations, that is
 \begin{equation*}
 S(p\circ T)=S(p),
 \end{equation*}
where $ T=\frac{ax+b}{cx+d}$ and $ad-bc\neq 0$. When
the Ricci tensor is parallel with respect to any of Berwald, Chern or Cartan
connection, it is constant along the geodesics and we can easily solve the equation (\ref{sch,Ric}), see \cite{BiS}.

Since $(M,F)$ is an Einstein space, plugging \eqref{Ric} in \eqref{sch,Ric} we get
\begin{equation}\label{S,K}
  S(p)=\frac{2}{n-1}F^2Ric=2F^2K,
 \end{equation}
 where $ K $ is a function of $ x $ alone.
If $ S(p)=0 $, then by the above equation, $K=0$, hence $ Ric=0 $, therefore by Proposition\ref{p1}, $(M,F)$ is a locally Minkowskian space.
 If $S(p)<0$, then  \eqref{S,K} yields $K<0$, and $Ric<0$, hence by Proposition\ref{p1}, $(M,F)$ is Riemannian.
\hspace {\stretch{1}}$\Box$

\textbf{ Acknowledgement.}
The first author would like to thank the Institut de Math\'{e}matiques de Toulouse (IMT) in which this article is partially written.

   Behroz Bidabad\\
Department of Mathematics and Computer Sciences\\
Amirkabir University of Technology (Tehran Polytechnic),
424 Hafez Ave. 15914 Tehran, Iran.
E-mail: bidabad@aut.ac.ir\\
Institut de Math\'{e}matique de Toulouse\\
Universit\'{e} Paul Sabatier, 118 route de Narbonne - F-31062 Toulouse, France.\\
behroz.bidabad@math.univ-toulouse.fr\\
Faranak Seddighi\\
Faculty of Mathematics, Payame Noor University of Tehran, Tehran, Iran.\\
F-Seddighi@student.pnu.ac.ir

\end{document}